\documentclass[12pt]{article}
\usepackage{amssymb}
\usepackage{amsfonts}
\usepackage{amsmath}

\input{tcilatex}
\begin{document}

\begin{center}
\textbf{A characterization of the fuzzy fractals}

\textbf{generated by an orbital fuzzy iterated function system}

\bigskip

by\textit{\ Radu MICULESCU,\ Alexandru MIHAIL }and\textit{\ Irina SAVU}

\bigskip

\textit{Dedicated to Professor Mihail MEGAN on the occasion of his 75}$^{%
\text{th}}$\textit{\ birthday}

\bigskip
\end{center}

\textbf{Abstract}. {\small Orbital fuzzy iterated function systems are
obtained as a combination of the concepts of iterated fuzzy set system and
orbital iterated function system. It turns out that, for such a system, the
corresponding fuzzy operator is weakly Picard, its fixed points being called
fuzzy fractals. In this paper we present a structure result concerning fuzzy
fractals associated to an orbital fuzzy iterated function system by proving
that such an object is perfectly determined by the action of the initial
term of the Picard iteration sequence on the closure of the orbits of
certain elements.}

\bigskip



\bigskip

\textbf{1. Introduction}

\bigskip

Fuzzy sets have their origin in Zadeh's remark that "more often than not,
the classes of objects encountered in the real physical world do not have
precisely defined criteria of membership" (see [16]). They have been
introduced, in 1965, with the purpose of reconciling mathematical modeling
and human knowledge in engineering sciences. To be precise, Zadeh was focus
on their potential applications "in human thinking, particularly in the
domain of pattern recognition, communication of information, and
abstraction".

Theory of iterated function systems, which was initiated, in 1981, by J.
Hutchinson (see [8]) and enriched by M. Barnsley and S. Demko (see [2] and
[3]), has at its core the construction of deterministic fractals and
measures. It has applications in image processing, stochastic growth model,
random dynamical systems, bioinformatics, economics, finance, engineering
sciences, human anatomy, physics etc.

The fuzzification Zadeh's idea was naturally adjusted to the
Hutchinson-Barnsley theory of iterated function systems. More precisely,
in1991, C. Cabrelli, B. Forte, U. Molter and E. Vrscay (see [4] and [5])
introduced the concept of iterated fuzzy set system which consists of a
finite family $(f_{i})_{i\in I}$ of contractions on a compact metric space $%
(X,d)$ together with a family of "grey level" maps $(\phi _{i})_{i\in I}$,
where $\phi _{i}:[0,1]\rightarrow \lbrack 0,1]$. One can associate to such a
system an operator on the class of normalized uppersemicontinuous fuzzy sets
of $X$ which turns out to be a contraction with respect to a metric $%
d_{\infty }$ involving the Hausdorff-Pompeiu distances between level sets.
Its unique fixed point is called the invariant fuzzy set. The relevance of
this theory to image processing is mentioned in [5]. The continuity
properties of the invariant fuzzy set with respect to changes in the
contractions $f_{i}$ and grey level maps $\phi _{i}$ are studied in [7].
Other papers dealing with iterated fuzzy set systems are [1] \ and [13]. Let
us also mention that R. Uthayakumar and D. Easwaramoorthy (see [15]) studied
the Hutchinson-Barnsley theory in the framework of the fuzzy hyperspace with
respect to the Hausdorrf-Pompeiu fuzzy metric.

Another natural generalization of Hutchinson' concept of iterated function
system, termed as orbital iterated function system, was recently considered
in [9], [10] and [14]. Here the idea is to consider iterated function
systems consisting of continuous functions satisfying Banach's orbital
condition. The novelty of this approach is that the associated fractal
operator is weakly Picard. It comes to the light that this is a genuine
generalization since there exist such systems for which the fractal operator
is weakly Picard, but not Picard (see Remark 4.1 from [9]). For extra
properties of this kind of system see [10] and [14].

A natural continuation of research lines previously mentioned is to examine
the so called orbital fuzzy iterated function systems (see [12]) which are
obtained as a combination of iterated fuzzy set systems and orbital iterated
function systems. It was proved (see Theorem 3.1 from [12]) that the
corresponding fuzzy operator is weakly Picard. Its fixed points are called
fuzzy fractals. More precisely, let us suppose that $\mathcal{S}%
_{Z}=((X,d),(f_{i})_{i\in I},(\rho _{i})_{i\in I})$ is such a system, $Z:%
\mathcal{F}_{\mathcal{S}}^{\ast }\rightarrow \mathcal{F}_{\mathcal{S}}^{\ast
}$ (where $\mathcal{F}_{\mathcal{S}}^{\ast }$ is a certain class of fuzzy
sets - see Section 2 for details) is the fuzzy Hutchinson-Barnsley operator
associated to $\mathcal{S}_{Z}$ and let us arbitrarily choose an element $u$
from $\mathcal{F}_{\mathcal{S}}^{\ast }$. Then the sequence $%
(Z^{[n]}(u))_{n\in \mathbb{N}}$ is convergent and its limit, denoted by $%
\mathbf{u}_{u}$, is a fuzzy fractal.

The goal of the present paper is to provide, for each $u\in \mathcal{F}_{%
\mathcal{S}}^{\ast }$, a description of the fuzzy fractal $\mathbf{u}_{u}$
in terms of certain fuzzy fractals $\mathbf{u}_{x}$ obtained as the limit of
the Picard iteration sequence which starts with a fuzzy set $u^{x}$
associated to $u$ and $x\in X$ such that $u(x)>0$. More precisely, Theorem
3.9, which is our main result, states that $\mathbf{u}_{u}=\underset{x\text{
such that }u(x)>0}{\max }\mathbf{u}_{x}=\underset{x\text{ such that }u(x)=1}{%
\max }\mathbf{u}_{x}$.

\bigskip

\textbf{2.} \textbf{Preliminaries}

\bigskip

\textbf{A}. \textbf{Basic} \textbf{notations and terminology}

\bigskip

By $\mathbb{N}$ we mean the set $\{1,2,...\}$.

For a family of functions $(f_{i})_{i\in I}$, where $f_{i}:X\rightarrow 
\mathbb{R}$, we shall use the following notation:

\begin{equation*}
\underset{i\in I}{\sup }f_{i}\overset{not}{=}\underset{i\in I}{\vee }f_{i}%
\text{.}
\end{equation*}

For a function $f:X\rightarrow X$ and $n\in \mathbb{N}$, the composition of $%
f$ by itself $n$ times is denoted by $f^{[n]}$.

A function $f:X\rightarrow X$, where $(X,d)$\ is a metric space, is called
weakly Picard operator if the sequence $(f^{[n]}(x))_{n\in \mathbb{N}}$ is
convergent for every $x\in X$\ and the limit (which may depend on $x$) is a
fixed point of $f$. A weakly Picard operator having a unique fixed point is
called Picard operator.

For a subset $A$ of a metric space $(X,d)$, by $diam(A)$ we mean the
diameter of $A$ i.e. $\underset{x,y\in A}{\sup }d(x,y)$.

For a metric space $(X,d)$ we shall use the following notations:

\begin{equation*}
\{A\subseteq X\mid A\neq \emptyset \text{ and }A\text{ is bounded}\}\overset{%
not}{=}P_{b}(X)
\end{equation*}

\begin{equation*}
\{A\subseteq X\mid A\neq \emptyset \text{ and }A\text{ is closed}\}\overset{%
not}{=}P_{cl}(X)
\end{equation*}

\begin{equation*}
P_{b}(X)\cap P_{cl}(X)\overset{not}{=}P_{b,cl}(X)
\end{equation*}%
\begin{equation*}
\{A\subseteq X\mid A\neq \emptyset \text{ and }A\text{ is compact}\}\overset{%
not}{=}P_{cp}(X)\text{.}
\end{equation*}

For a metric space $(X,d)$, by $h$ we designate the Hausdorff-Pompeiu metric
on $X$, i.e. the function $h:P_{b,cl}(X)\times P_{b,cl}(X)\rightarrow
\lbrack 0,\infty )$, described by 
\begin{equation*}
h(K_{1},K_{2})=\max \{\underset{x\in K_{1}}{\sup }d(x,K_{2}),\underset{x\in
K_{2}}{\sup }d(x,K_{1})\}\text{,}
\end{equation*}%
for every $K_{1},K_{2}\in P_{b,cl}(X)$.

\bigskip

\textbf{Remark 2.1.} $(P_{cp}(X),h)$ \textit{is a complete metric space,
provided that }$(X,d)$\textit{\ is complete. If }$(A_{n})_{n\in \mathbb{N}%
}\subseteq P_{b,cl}(X)$\textit{\ is Cauchy, then }$\underset{n\rightarrow
\infty }{\lim }A_{n}=\{x\in X\mid $\textit{there exists a strictly
increasing sequence }$(n_{k})_{k\in \mathbb{N}}\subseteq \mathbb{N}$\textit{%
\ and }$x_{n_{k}}\in A_{n_{k}}$\textit{\ for every }$k\in \mathbb{N}$\textit{%
\ such that }$\underset{k\rightarrow \infty }{\lim }x_{n_{k}}=x\}$\textit{.}

\bigskip

\textbf{B}. \textbf{Fuzzy sets}

\bigskip

For a set $X$, we shall use the following notation:%
\begin{equation*}
\{u:X\rightarrow \lbrack 0,1]\}\overset{not}{=}\mathcal{F}_{X}
\end{equation*}

The elements of $\mathcal{F}_{X}$ are called fuzzy subsets of $X$.

A non-zero function $\rho :[0,1]\rightarrow \lbrack 0,1]$ is called a grey
level map.

To every grey level map $\rho $ and $u\in \mathcal{F}_{X}$ one could
associate the element of $\mathcal{F}_{X}$, denoted by $\rho (u)$, given by $%
\rho \circ u$.

$u\in \mathcal{F}_{X}$ is called normal if there exists $x\in X$ such that $%
u(x)=1$.

For $u\in \mathcal{F}_{X}$ and $\alpha \in (0,1]$, we shall use the
following notations:%
\begin{equation*}
\{x\in X\mid u(x)\geq \alpha \}\overset{not}{=}[u]^{\alpha }
\end{equation*}%
\begin{equation*}
\{x\in X\mid u(x)>0\}\overset{not}{=}[u]^{\ast }\text{.}
\end{equation*}

Given a metric space $(X,d)$, $u\in \mathcal{F}_{X}$ is called compactly
supported if supp$u:=\overline{[u]^{\ast }}\overset{not}{=}[u]^{0}\in
P_{cp}(X)$.

For a metric space $(X,d)$, we shall use the following notations:%
\begin{equation*}
\{u\in \mathcal{F}_{X}\mid u\text{ is normal and compactly supported}\}%
\overset{not}{=}\mathcal{F}_{X}^{\ast \ast }
\end{equation*}%
\begin{equation*}
\{u\in \mathcal{F}_{X}^{\ast \ast }\mid u\text{ is upper semicontinuous}\}%
\overset{not}{=}\mathcal{F}_{X}^{\ast }\text{.}
\end{equation*}

For a metric space $(X,d)$ and $x\in X$ we consider $\delta _{x}\in \mathcal{%
F}_{X}^{\ast }$ given by 
\begin{equation*}
\delta _{x}(t)=\{%
\begin{array}{cc}
1\text{,} & \text{if }t=x \\ 
0\text{,} & \text{if }t\neq x%
\end{array}%
\text{.}
\end{equation*}

\bigskip

\textbf{Remark 2.2.}%
\begin{equation*}
\text{supp}\delta _{x}=\{x\}\text{.}
\end{equation*}

\bigskip

To every $f:X\rightarrow Y$ and $u\in \mathcal{F}_{X}$ one could associate
an element of $\mathcal{F}_{Y}$, denoted by $f(u)$, which is described in
the following way:%
\begin{equation*}
f(u)(y)=\{%
\begin{array}{cc}
\underset{x\in f^{-1}(\{y\})}{\sup }u(x)\text{,} & \text{if }%
f^{-1}(\{y\})\neq \emptyset \\ 
0\text{,} & \text{if }f^{-1}(\{y\})=\emptyset%
\end{array}%
\text{,}
\end{equation*}%
for every $y\in Y$.

For a metric space $(X,d)$, the function $d_{\infty }:\mathcal{F}_{X}^{\ast
\ast }\times \mathcal{F}_{X}^{\ast \ast }\rightarrow \lbrack 0,\infty ]$,
given by%
\begin{equation*}
d_{\infty }(u,v)\overset{def}{=}\underset{\alpha \in \lbrack 0,1]}{\sup }%
h([u]^{\alpha },[v]^{\alpha })\overset{\text{Lemma 2.5 from [12] }}{=}%
\underset{\alpha \in (0,1]}{\sup }h([u]^{\alpha },[v]^{\alpha })\text{,}
\end{equation*}%
for every $u,v\in \mathcal{F}_{X}^{\ast \ast }$, is semidistance on $%
\mathcal{F}_{X}^{\ast \ast }$. Its restriction to $\mathcal{F}_{X}^{\ast
}\times \mathcal{F}_{X}^{\ast }$ is a metric on $\mathcal{F}_{X}^{\ast }$
(see [6]), which, for the sake of simplicity, will be also denoted by $%
d_{\infty }$. Moreover $(\mathcal{F}_{X}^{\ast },d_{\infty })$ is a complete
metric space provided that the metric space $(X,d)$ is complete.

\bigskip

\textbf{C}. \textbf{Iterated function systems}

\bigskip

An iterated function system (IFS for short)\ consists of:

i) a complete metric space $(X,d)$;

ii) a finite family of contractions $f_{i}:X\rightarrow X$, with $i\in I$.

We denote by $\mathcal{S}=((X,d),(f_{i})_{i\in I})$ such an IFS.

One can associate to such a system $\mathcal{S}=((X,d),(f_{i})_{i\in I})$
the function $F_{\mathcal{S}}:P_{cp}(X)\rightarrow P_{cp}(X)$, given\ by

\begin{equation*}
F_{\mathcal{S}}(K)=\underset{i\in I}{\cup }f_{i}(K)\text{,}
\end{equation*}%
for all $K\in P_{cp}(X)$, which is called the fractal operator associated to 
$\mathcal{S}$.

It turns out (see [8]) that $F_{\mathcal{S}}$ is a Banach contraction on the
complete metric space $(P_{cp}(X),h)$, so it is a Picard operator with
respect to $h$ and its fixed point (which is denoted by $A_{\mathcal{S}}$)\
is called the attractor of\ $\mathcal{S}$.

\bigskip

\textbf{D.} \textbf{Iterated fuzzy function systems}

\bigskip

An iterated fuzzy function system\ consists of:

i) an iterated function system $\mathcal{S}=((X,d),(f_{i})_{i\in I})$

ii) an admissible system of grey level maps $(\rho _{i})_{i\in I}$ i.e. $%
\rho _{i}(0)=0$, $\rho _{i}$ is nondecreasing and right continuous for every 
$i\in I$ and there exists $j\in I$ such that $\rho _{j}(1)=1$.

We denote by $\mathcal{S}_{Z}=((X,d),(f_{i})_{i\in I},(\rho _{i})_{i\in I})$
such a system.

One can associate to such a system $\mathcal{S}_{Z}=((X,d),(f_{i})_{i\in
I},(\rho _{i})_{i\in I})$ the function $Z:\mathcal{F}_{X}^{\ast }\rightarrow 
\mathcal{F}_{X}^{\ast }$, given\ by

\begin{equation*}
Z(u)=\underset{i\in I}{\vee }\rho _{i}(f_{i}(u))\text{,}
\end{equation*}%
for all $u\in \mathcal{F}_{X}^{\ast }$, which is called the fuzzy
Hutchinson-Barnsley operator associated to $\mathcal{S}_{Z}$. Note that $Z$\
is well defined (see Proposition 2.12 from [13]).

It turns out (see Theorem 2.14 from [13]) that $Z$ is a Banach contraction
on the complete metric space $(\mathcal{F}_{X}^{\ast },d_{\infty })$ (so it
is a Picard operator) whose unique fixed point is called the fuzzy fractal
generated by $\mathcal{S}_{Z}$ (note that its support is a subset of $A_{%
\mathcal{S}}$ - see Theorem 2.4.2 from [5] or Theorem 2.21 from [13]).

\bigskip

\textbf{E.} \textbf{Orbital iterated function systems}

\bigskip

An orbital iterated function system\ consists of:

i) a complete metric space $(X,d)$;

ii) a finite family of continuous functions $f_{i}:X\rightarrow X$, $i\in I$%
, having the property that there exists $C\in \lbrack 0,1)$ such that $%
d(f_{i}(y),f_{i}(z))\leq Cd(y,z)$ for every $i\in I$, $x\in X$ and $y,z\in 
\mathcal{O}(x)$, where, for $B\in P_{cp}(X)$, by the orbit of $B$, denoted
by $\mathcal{O}(B)$, we mean the set $B\cup \underset{n\in \mathbb{N}\text{, 
}\omega _{1},..,\omega _{1}\in I}{\cup }(f_{\omega _{1}}\circ f_{\omega
_{2}}\circ ...\circ f_{\omega _{n}})(B)$ and we adopt the notation $\mathcal{%
O}(\{x\})\overset{not}{=}\mathcal{O}(x)$ for every $x\in X$.

We denote by $\mathcal{S}=((X,d),(f_{i})_{i\in I})$ such a system.

As for the case of IFSs, one can associate to an orbital iterated function
system $\mathcal{S}$ its fractal operator. It turns out (see [9] and [11])
that the fractal operator associated to an orbital function system is a
weakly Picard operator with respect to the Hausdorff-Pompeiu metric, every
of its fixed points being called an attractor of the system.

If $\mathcal{S}=((X,d),(f_{i})_{i\in I})$ is an orbital iterated function
system, $K\in P_{cp}(X)$ and $x\in X$, we shall use the following notations:

\begin{equation*}
\underset{n\rightarrow \infty }{\lim }F_{S}^{[n]}(K)\overset{not}{=}A_{K}
\end{equation*}%
and%
\begin{equation*}
A_{\{x\}}\overset{not}{=}A_{x}\text{.}
\end{equation*}

\bigskip

\textbf{F.} \textbf{Orbital fuzzy iterated function systems}

\bigskip

An orbital fuzzy iterated function system\ consists of:

i) an orbital iterated function system $((X,d),(f_{i})_{i\in I})$

ii) an admissible system of grey level maps $(\rho _{i})_{i\in I}$ i.e. $%
\rho _{i}(0)=0$, $\rho _{i}$ is nondecreasing and right continuous for every 
$i\in I$ and there exists $j\in I$ such that $\rho _{j}(1)=1$.

We denote by $\mathcal{S}_{Z}=((X,d),(f_{i})_{i\in I},(\rho _{i})_{i\in I})$
such a system.

One can associate to such a system $\mathcal{S}_{Z}=((X,d),(f_{i})_{i\in
I},(\rho _{i})_{i\in I})$ the function $Z:\mathcal{F}_{X}^{\ast \ast
}\rightarrow \mathcal{F}_{X}^{\ast \ast }$, given\ by

\begin{equation*}
Z(u)=\underset{i\in I}{\vee }\rho _{i}(f_{i}(u))\text{,}
\end{equation*}%
for all $u\in \mathcal{F}_{X}^{\ast \ast }$, which is called the fuzzy
Hutchinson-Barnsley operator associated to $\mathcal{S}_{Z}$. Note that, in
view of Proposition 2.12 from [13], $Z$\ is well defined.

We shall use the following notations:%
\begin{equation*}
\{u\in \mathcal{F}_{X}^{\ast \ast }\mid \text{for each }x\in \lbrack
u]^{\ast }\text{ there exist}
\end{equation*}%
\begin{equation*}
w_{x},y_{x}\in X\text{ such that }x,y_{x}\in \mathcal{O}(w_{x})\text{ and }%
u(y_{x})=1\}\overset{not}{=}\mathcal{F}_{\mathcal{S}}^{\ast \ast }
\end{equation*}%
and%
\begin{equation*}
\{u\in \mathcal{F}_{\mathcal{S}}^{\ast \ast }\mid u\text{ is upper
semicontinuous}\}\overset{not}{=}\mathcal{F}_{\mathcal{S}}^{\ast }\text{.}
\end{equation*}

Note that 
\begin{equation*}
\delta _{x}\in \mathcal{F}_{\mathcal{S}}^{\ast }\text{,}
\end{equation*}
for every $x\in X$.

For $u\in \mathcal{F}_{\mathcal{S}}^{\ast }$ and $x\in \lbrack u]^{\ast }$
(hence there exist $w_{x},y_{x}\in X$ such that $x,y_{x}\in \mathcal{O}%
(w_{x})$ and $u(y_{x})=1$), we shall use the following notations:%
\begin{equation*}
\underset{n\rightarrow \infty }{\lim }Z^{[n]}(u)\overset{not}{=}\mathbf{u}%
_{u}\overset{\text{Lemma 3.1 from [12] \& Remark 2.3}}{\in }\mathcal{F}_{%
\mathcal{S}}^{\ast }\text{,}
\end{equation*}%
\begin{equation*}
\underset{n\rightarrow \infty }{\lim }Z^{[n]}(u^{x})\overset{not}{=}\mathbf{u%
}_{x}\overset{\text{Lemma 3.1 from [12] \& Remark 2.3}}{\in }\mathcal{F}_{%
\mathcal{S}}^{\ast }
\end{equation*}%
and%
\begin{equation*}
w_{u}\overset{not}{=}\underset{x\in \lbrack u]^{\ast }}{\vee }\mathbf{u}_{x}%
\text{,}
\end{equation*}%
where $u^{x}\in \mathcal{F}_{\mathcal{S}}^{\ast }$ is described by%
\begin{equation*}
u^{x}(y)=\{%
\begin{array}{cc}
u(y)\text{,} & \text{if }y\in \overline{\mathcal{O}(w_{x})} \\ 
0\text{,} & \text{otherwise}%
\end{array}%
\text{.}
\end{equation*}%
Note that the existence of the above limits is based on Remark 2.5 and that,
according to Proposition 3.1, $\mathbf{u}_{x}$ is well defined.

\bigskip

\textbf{Remark 2.3 }(see Proposition 2.11 and Lemma 3.3 from [12])\textbf{.} 
\textit{In the above framework,} \textit{we have}

\begin{equation*}
Z(\mathcal{F}_{X}^{\ast })\subseteq \mathcal{F}_{X}^{\ast }\text{, }Z(%
\mathcal{F}_{\mathcal{S}}^{\ast \ast })\subseteq \mathcal{F}_{\mathcal{S}%
}^{\ast \ast }\text{ \textit{and} }Z(\mathcal{F}_{\mathcal{S}}^{\ast
})\subseteq \mathcal{F}_{\mathcal{S}}^{\ast }\text{.}
\end{equation*}

\bigskip

\textbf{Remark 2.4}\textit{\ }(see Claim 3.5 from the proof of Theorem 3.1
from [12])\textbf{.} \textit{In the above framework,\ }$\mathbf{Z}:\mathcal{F%
}_{\mathcal{S}}^{\ast }\rightarrow \mathcal{F}_{\mathcal{S}}^{\ast }$\textit{%
, given\ by} $\mathbf{Z}(u)=Z(u)$ \textit{for every} $u\in \mathcal{F}_{%
\mathcal{S}}^{\ast }$,\textit{\ is continuous.}

\bigskip

\textbf{Remark 2.5}\textit{\ }(see Theorem 3.1 from [12])\textbf{.} \textit{%
In the above framework,} $\mathbf{Z}$\textit{\ is weakly Picard. Its fixed
points are called fuzzy fractals generated by the orbital fuzzy iterated
function system }$\mathcal{S}_{Z}$\textit{.}

\bigskip

\textbf{Remark 2.6}\textit{\ }(see Lemma 3.4 from [12])\textbf{. }\textit{In
the above framework,} \textit{for each family }$(u_{j})_{j\in J}$\textit{\
of elements from }$\mathcal{F}_{X}^{\ast \ast }$\textit{, where }$J$ \textit{%
is infinite, we have}%
\begin{equation*}
\underset{j\in J}{\vee }Z(u_{j})=\underset{j\in J}{\max }\text{ }Z(u_{j})%
\text{,}
\end{equation*}%
\textit{provided that:}

\textit{i) there exists }$K\in P_{cp}(X)$ \textit{such that }supp$%
u_{j}\subseteq K$\textit{\ for all }$j\in J$\textit{;}

\textit{ii) }$\underset{j\in J}{\vee }u_{j}=\underset{j\in J}{\max }$ $u_{j}$%
;

\textit{iii) }$\underset{j\in J}{\vee }u_{j}\in \mathcal{F}_{X}^{\ast }$.

\bigskip

\textbf{Remark 2.7}\textit{\ }(see Lemma 3.5 from [12])\textbf{.} \textit{In
the above framework, for every families} $(u_{j})_{j\in J}$\textit{\ and }$%
(v_{j})_{j\in J}$ \textit{of elements from }$\mathcal{F}_{X}^{\ast \ast }$%
\textit{, where }$J$ \textit{is infinite, we have}%
\begin{equation*}
d_{\infty }(\underset{j\in J}{\vee }u_{j},\underset{j\in J}{\vee }v_{j})\leq 
\underset{j\in J}{\sup }\text{ }d_{\infty }(u_{j},v_{j})\text{,}
\end{equation*}%
\textit{provided that:}

\textit{i) there exists }$K\in P_{cp}(X)$ \textit{such that }supp$%
u_{j}\subseteq K$ \textit{and} supp$v_{j}\subseteq K$\textit{\ for all }$%
j\in J$\textit{;}

\textit{ii) }$\underset{j\in J}{\vee }u_{j}=\underset{j\in J}{\max }$ $u_{j}$
and $\underset{j\in J}{\vee }v_{j}=\underset{j\in J}{\max }$ $v_{j}$.

\bigskip

\textbf{Remark 2.8. }\textit{In the above framework,}%
\begin{equation*}
\text{supp}\mathbf{u}_{x}\subseteq A_{x}\text{,}
\end{equation*}%
\textit{for every} $u\in \mathcal{F}_{\mathcal{S}}^{\ast }$ \textit{and} $%
x\in \lbrack u]^{\ast }$.

Indeed, it results from the following relations:%
\begin{equation*}
\text{supp}\mathbf{u}_{x}\overset{\text{Remark 2.4}}{=}\text{supp}Z(\mathbf{u%
}_{x})\overset{\text{Lemma 3.7 from [12]}}{\subseteq }F_{\mathcal{S}}(\text{%
supp}\mathbf{u}_{x})\text{.}
\end{equation*}

It also can be derived from Theorem 2.4.2 from [5].

\bigskip

\textbf{Remark 2.9} (see Lemma 4.5 from [11)\textbf{.} \textit{In the above
framework, for every }$x_{1},x_{2}\in X$\textit{, we have}%
\begin{equation*}
A_{x_{1}}=A_{x_{2}}
\end{equation*}%
\textit{provided that} $\mathcal{O}(x_{1})\cap \mathcal{O}(x_{1})\neq
\emptyset $\textit{.}

\bigskip

\textbf{Remark 2.10}\textit{\ }(see Remark 2.1 from [12])\textbf{.} \textit{%
In the above framework,} 
\begin{equation*}
\overline{\mathcal{O}(x)}=\mathcal{O}(x)\cup A_{x}\text{,}
\end{equation*}%
\textit{for every} $x\in X$\textit{.}

\bigskip

\textbf{Remark 2.11.} \textit{In the above framework,}%
\begin{equation*}
w_{u}\in \mathcal{F}_{X}^{\ast \ast }\text{,}
\end{equation*}%
\textit{for every} $u\in \mathcal{F}_{\mathcal{S}}^{\ast }$.

Indeed, on the one hand, as $\mathbf{u}_{x}\in \mathcal{F}_{X}^{\ast }$ is
normal for every $x\in \lbrack u]^{\ast }$, we deduce that $w_{u}$ is
normal. On the other hand, since supp$\mathbf{u}_{x}\overset{\text{Remark 2.8%
}}{\subseteq }A_{x}\overset{\text{Proposition 5 from [14]}}{\subseteq }A_{%
\text{supp}u}\in P_{cp}(X)$ for every $x\in \lbrack u]^{\ast }$, we conclude
that supp$w_{u}=$supp$\underset{x\in \lbrack u]^{\ast }}{\vee }\mathbf{u}%
_{x}\subseteq A_{\text{supp}u}\in P_{cp}(X)$, so supp$w_{u}\in P_{cp}(X)$.

\bigskip

\textbf{Remark 2.12.} \textit{In the above framework,}%
\begin{equation*}
\text{supp}u^{x}\subseteq \text{supp}u\text{,}
\end{equation*}%
\textit{for every} $u\in \mathcal{F}_{\mathcal{S}}^{\ast }$ \textit{and }$%
x\in \lbrack u]^{\ast }$.

\bigskip

\textbf{Remark 2.13 }(see Claim 3.2 from the proof of Theorem 3.1 from [12])%
\textbf{.} \textit{In the above framework, we have:}

\textit{a)}%
\begin{equation*}
Z^{[n]}(u)=\underset{x\in \lbrack u]^{\ast }}{\vee }Z^{[n]}(u^{x})\text{,}
\end{equation*}%
\textit{for every }$n\in \mathbb{N}$ \textit{and every }$u\in \mathcal{F}_{%
\mathcal{S}}^{\ast }$\textit{.}

\textit{b)}%
\begin{equation*}
d_{\infty }(Z^{[n]}(u),\mathbf{u}_{u})\leq \frac{C^{n}}{1-C}diam(F_{S}(\text{%
supp}u)\cup \text{supp}u)\text{,}
\end{equation*}%
\textit{for all }$n\in \mathbb{N}$ \textit{and} $u\in \mathcal{F}_{\mathcal{S%
}}^{\ast }$\textit{.}

\bigskip

\textbf{3. The main results}

\bigskip

\textbf{Proposition 3.1.} \textit{Let} $\mathcal{S}_{Z}=((X,d),(f_{i})_{i\in
I},(\rho _{i})_{i\in I})$ \textit{be an orbital fuzzy iterated function
system, }$u\in \mathcal{F}_{\mathcal{S}}^{\ast }$ \textit{and} $x\in \lbrack
u]^{\ast }$ \textit{(hence there exist }$w_{x},y_{x}\in X$\textit{\ such
that }$x,y_{x}\in \mathcal{O}(w_{x})$\textit{\ and }$u(y_{x})=1$\textit{).
Then }%
\begin{equation*}
\underset{n\rightarrow \infty }{\lim }Z^{[n]}(\delta _{s})=\mathbf{u}_{x}%
\text{,}
\end{equation*}%
\textit{for every} $s\in \overline{\mathcal{O}(w_{x})}$.\textit{\ In
particular} 
\begin{equation*}
\underset{n\rightarrow \infty }{\lim }Z^{[n]}(\delta _{x})=\mathbf{u}_{x}%
\text{.}
\end{equation*}

\textit{Proof}. Because $((\overline{\mathcal{O}(w_{x})},d),(\overset{\sim }{%
f_{i}})_{i\in I},(\rho _{i})_{i\in I})$, where $\overset{\sim }{f_{i}}:%
\overline{\mathcal{O}(w_{x})}\rightarrow \overline{\mathcal{O}(w_{x})}$ is
given by $\overset{\sim }{f_{i}}(y)=f_{i}(y)$ for every $y\in \overline{%
\mathcal{O}(w_{x})}$, has a unique fuzzy fractal (as it is an iterated fuzzy
function system) and $\delta _{s\mid X\smallsetminus \overline{\mathcal{O}%
(w_{x})}}=u_{\mid X\smallsetminus \overline{\mathcal{O}(w_{x})}%
}^{x}=Z^{[n]}(\delta _{s\mid X\smallsetminus \overline{\mathcal{O}(w_{x})}%
})=Z^{[n]}(u_{\mid X\smallsetminus \overline{\mathcal{O}(w_{x})}}^{x})=0$
for every $s\in \overline{\mathcal{O}(w_{x})}$ and every $n\in \mathbb{N}$,
we deduce that%
\begin{equation*}
\underset{n\rightarrow \infty }{\lim }Z^{[n]}(\delta _{s})=\underset{%
n\rightarrow \infty }{\lim }Z^{[n]}(u^{x})=\mathbf{u}_{x}\text{,}
\end{equation*}%
for every $s\in \overline{\mathcal{O}(w_{x})}$. $\square $

\bigskip

\textbf{Proposition 3.2.} \textit{Let} $\mathcal{S}_{Z}=((X,d),(f_{i})_{i\in
I},(\rho _{i})_{i\in I})$ \textit{be an orbital fuzzy iterated function
system} \textit{and} $u\in \mathcal{F}_{\mathcal{S}}^{\ast }$\textit{. Then }%
\begin{equation*}
\mathbf{u}_{y}=\mathbf{u}_{x}\text{,}
\end{equation*}%
\textit{for every} $x\in \lbrack u]^{\ast }$\ \textit{and every} $y\in
\lbrack \mathbf{u}_{x}]^{\ast }$.

\textit{Proof}. Let us consider $x\in \lbrack u]^{\ast }$\ and $y\in \lbrack 
\mathbf{u}_{x}]^{\ast }$. As $u\in \mathcal{F}_{\mathcal{S}}^{\ast }$ ,
there exist $w_{x},y_{x}\in X$ such that 
\begin{equation}
x,y_{x}\in \mathcal{O}(w_{x})  \tag{1}
\end{equation}%
and $u(y_{x})=1$.

In addition 
\begin{equation}
y\in supp\mathbf{u}_{x}\overset{\text{Remark 2.8}}{\subseteq }A_{x}\overset{%
\text{(1) \& Remark 2.9}}{=}A_{w_{x}}\overset{\text{Remark 2.10}}{\subseteq }%
\overline{\mathcal{O}(w_{x})}\text{.}  \tag{2}
\end{equation}

Therefore%
\begin{equation*}
\mathbf{u}_{x}\overset{\text{(2) \& Proposition 3.1}}{=}\underset{%
n\rightarrow \infty }{\lim }Z^{[n]}(\delta _{y})\overset{\text{Proposition
3.1}}{=}\mathbf{u}_{y}\text{. }\square
\end{equation*}

\bigskip

\textbf{Proposition 3.3.} \textit{Let} $\mathcal{S}_{Z}=((X,d),(f_{i})_{i\in
I},(\rho _{i})_{i\in I})$ \textit{be an orbital fuzzy iterated function
system} \textit{and} $u\in \mathcal{F}_{\mathcal{S}}^{\ast }$\textit{. Then
the function} $U:[u]^{\ast }\rightarrow \mathcal{F}_{X}^{\ast }$, \textit{%
given by} 
\begin{equation*}
U(x)=\mathbf{u}_{x}\text{,}
\end{equation*}%
\textit{for every} $x\in \lbrack u]^{\ast }$, \textit{is continuous.}

\textit{Proof}. We are going to prove that $U$ is sequentially continuous.

We consider $(x_{n})_{n\in \mathbb{N}}\subseteq \lbrack u]^{\ast }$ and $%
x\in \lbrack u]^{\ast }$ such that $\underset{n\rightarrow \infty }{\lim }%
x_{n}=x$ and we will prove that $\underset{n\rightarrow \infty }{\lim }%
U(x_{n})=U(x)$, i.e. $\underset{n\rightarrow \infty }{\lim }\mathbf{u}%
_{x_{n}}=\mathbf{u}_{x}$.

We have%
\begin{equation}
\text{supp}\delta _{x}\overset{\text{Remark 2.2}}{=}\{x\}\subseteq K\text{
and supp}\delta _{x_{n}}\overset{\text{Remark 2.2}}{=}\{x_{n}\}\subseteq K%
\text{,}  \tag{1}
\end{equation}%
for every $n\in \mathbb{N}$, where 
\begin{equation*}
\{x_{n}\mid n\in \mathbb{N}\}\cup \{x\}\overset{not}{=}K\in P_{cp}(X)\text{.}
\end{equation*}

Since 
\begin{equation*}
\underset{n\rightarrow \infty }{\lim }d_{\infty }(\delta _{x_{n}},\delta
_{x})=\underset{n\rightarrow \infty }{\lim }\underset{\alpha \in (0,1]}{\sup 
}h([\delta _{x_{n}}]^{\alpha },[\delta _{x}]^{\alpha })=
\end{equation*}%
\begin{equation*}
=\underset{n\rightarrow \infty }{\lim }h(\{x_{n}\},\{x\})=\underset{%
n\rightarrow \infty }{\lim }d(x_{n},x)=0\text{,}
\end{equation*}%
via Remark 2.4, we infer that%
\begin{equation}
\underset{n\rightarrow \infty }{\lim }d_{\infty }(Z^{[m]}(\delta
_{x_{n}}),Z^{[m]}(\delta _{x}))=0\text{,}  \tag{2}
\end{equation}%
for every $m\in \mathbb{N}$.

The equality 
\begin{equation*}
\underset{n\rightarrow \infty }{\lim }Z^{[n]}(\delta _{s})\overset{\text{%
Proposition 3.1}}{=}\mathbf{u}_{s}\text{,}
\end{equation*}%
which is valid for every $s\in \lbrack u]^{\ast }$, leads to the conclusion
that%
\begin{equation}
d_{\infty }(\mathbf{u}_{s},Z^{[n]}(\delta _{s}))\overset{\text{Remark 2.2 \&
Remark 2.13, b)}}{\leq }\frac{C^{n}}{1-C}diam(F_{S}(\{s\})\cup \{s\})\text{,}
\tag{3}
\end{equation}%
for every $s\in \lbrack u]^{\ast }$ and every $n\in \mathbb{N}$.

Note that%
\begin{equation*}
d_{\infty }(\mathbf{u}_{x_{n}},\mathbf{u}_{x})\leq
\end{equation*}%
\begin{equation*}
\leq d_{\infty }(\mathbf{u}_{x_{n}},Z^{[m]}(\delta _{x_{n}}))+d_{\infty
}(Z^{[m]}(\delta _{x_{n}}),Z^{[m]}(\delta _{x}))+d_{\infty }(Z^{[m]}(\delta
_{x}),\mathbf{u}_{x})\overset{(3)}{\leq }
\end{equation*}%
\begin{equation}
\leq 2\frac{C^{m}}{1-C}diam(F_{S}(K)\cup K)+d_{\infty }(Z^{[m]}(\delta
_{x_{n}}),Z^{[m]}(\delta _{x}))\text{,}  \tag{4}
\end{equation}%
for every $m,n\in \mathbb{N}$.

Let us consider a fixed $\varepsilon >0$, but arbitrarily chosen.

As $\underset{m\rightarrow \infty }{\lim }2\frac{C^{m}}{1-C}%
diam(F_{S}(K)\cup K)=0$, there exists $m_{0}\in \mathbb{N}$ such that $2%
\frac{C^{m_{0}}}{1-C}diam(F_{S}(K)\cup K)<\frac{\varepsilon }{2}$ and, via $%
(4)$, we obtain%
\begin{equation}
d_{\infty }(\mathbf{u}_{x_{n}},\mathbf{u}_{x})\leq \frac{\varepsilon }{2}%
+d_{\infty }(Z^{[m_{0}]}(\delta _{x_{n}}),Z^{[m_{0}]}(\delta _{x}))\text{,} 
\tag{5}
\end{equation}%
for every $n\in \mathbb{N}$.

Since $\underset{n\rightarrow \infty }{\lim }d_{\infty }(Z^{[m_{0}]}(\delta
_{x_{n}}),Z^{[m_{0}]}(\delta _{x}))\overset{(2)}{=}0$, there exists $%
n_{\varepsilon }\in \mathbb{N}$ such that%
\begin{equation}
d_{\infty }(Z^{[m_{0}]}(\delta _{x_{n}}),Z^{[m_{0}]}(\delta _{x}))<\frac{%
\varepsilon }{2}\text{,}  \tag{6}
\end{equation}%
for every $n\in \mathbb{N}$, $n\geq n_{\varepsilon }$.

Using $(5)$ and $(6)$, we get%
\begin{equation*}
d_{\infty }(\mathbf{u}_{x_{n}},\mathbf{u}_{x})<\varepsilon \text{,}
\end{equation*}%
for every $n\in \mathbb{N}$, $n\geq n_{\varepsilon }$ which proves that $%
\underset{n\rightarrow \infty }{\lim }\mathbf{u}_{x_{n}}=\mathbf{u}_{x}$. $%
\square $

\bigskip

\textbf{Proposition 3.4.} \textit{Let} $\mathcal{S}_{Z}=((X,d),(f_{i})_{i\in
I},(\rho _{i})_{i\in I})$ \textit{be an orbital fuzzy iterated function
system and} $u\in \mathcal{F}_{\mathcal{S}}^{\ast }$\textit{. Then} 
\begin{equation*}
\lbrack \underset{x\in \lbrack u]^{\ast }}{\vee }\mathbf{u}_{x}]^{\alpha }=%
\underset{x\in \lbrack u]^{\ast }}{\cup }[\mathbf{u}_{x}]^{\alpha }\text{,}
\end{equation*}%
\textit{for every} $\alpha \in (0,1]$ \textit{and} $u\in \mathcal{F}_{%
\mathcal{S}}^{\ast }$\textit{.}

\textit{Proof}. Let us consider a fixed $\alpha \in (0,1]$, but arbitrarily
chosen.

First we prove the inclusion%
\begin{equation}
\lbrack \underset{x\in \lbrack u]^{\ast }}{\vee }\mathbf{u}_{x}]^{\alpha
}\subseteq \underset{x\in \lbrack u]^{\ast }}{\cup }[\mathbf{u}_{x}]^{\alpha
}\text{.}  \tag{1}
\end{equation}

For $y\in \lbrack \underset{x\in \lbrack u]^{\ast }}{\vee }\mathbf{u}%
_{x}]^{\alpha }$ we have%
\begin{equation}
\underset{x\in \lbrack u]^{\ast }}{\sup }\mathbf{u}_{x}(y)\geq \alpha 
\tag{2}
\end{equation}%
and consequently there exists $x\in \lbrack u]^{\ast }$ such that $\mathbf{u}%
_{x}(y)>0$, i.e. $y\in \lbrack \mathbf{u}_{x}]^{\ast }$. Moreover
Proposition 3.2 ensures us that $\mathbf{u}_{x}(y)=\mathbf{u}_{y}(y)$ for
every $x\in \lbrack u]^{\ast }$ such that $\mathbf{u}_{x}(y)>0$ and thus,
via $(2)$, we conclude that $\mathbf{u}_{x}(y)\geq \alpha $ (so $y\in
\lbrack \mathbf{u}_{x}]^{\alpha }$) for every $x\in \lbrack u]^{\ast }$ such
that $\mathbf{u}_{x}(y)>0$. Consequently $y\in \underset{x\in \lbrack
u]^{\ast }}{\cup }[\mathbf{u}_{x}]^{\alpha }$ and the justification of $(1)$
is finished.

Now we prove the inclusion%
\begin{equation}
\underset{x\in \lbrack u]^{\ast }}{\cup }[\mathbf{u}_{x}]^{\alpha }\subseteq
\lbrack \underset{x\in \lbrack u]^{\ast }}{\vee }\mathbf{u}_{x}]^{\alpha }%
\text{.}  \tag{3}
\end{equation}

For $y\in \underset{x\in \lbrack u]^{\ast }}{\cup }[\mathbf{u}_{x}]^{\alpha
} $ there exists $x_{y}\in \lbrack u]^{\ast }$ such that $y\in \lbrack 
\mathbf{u}_{x_{y}}]^{\alpha }$, i.e. $\mathbf{u}_{x_{y}}(y)\geq \alpha $.
Hence $\underset{x\in \lbrack u]^{\ast }}{\sup }\mathbf{u}_{x}(y)\geq 
\mathbf{u}_{x_{y}}(y)\geq \alpha $, i.e. $y\in \lbrack \underset{x\in
\lbrack u]^{\ast }}{\vee }\mathbf{u}_{x}]^{\alpha }$, and the justification
of $(3)$ is finalized.

In view of $(1)$ and $(3)$ the proof is completed. $\square $

\bigskip

As a by-product of the above Proposition's proof we have the following

\bigskip

\textbf{Remark 3.5. }\textit{Let} $\mathcal{S}_{Z}=((X,d),(f_{i})_{i\in
I},(\rho _{i})_{i\in I})$ \textit{be an orbital fuzzy iterated function
system. Then} 
\begin{equation*}
\underset{x\in \lbrack u]^{\ast }}{\vee }\mathbf{u}_{x}=\underset{x\in
\lbrack u]^{\ast }}{\max }\text{ }\mathbf{u}_{x}\text{,}
\end{equation*}%
\textit{for every} $u\in \mathcal{F}_{\mathcal{S}}^{\ast }$\textit{.}

\bigskip

\textbf{Proposition 3.6.} \textit{Let} $\mathcal{S}_{Z}=((X,d),(f_{i})_{i\in
I},(\rho _{i})_{i\in I})$ \textit{be an orbital fuzzy iterated function
system and }$u\in \mathcal{F}_{\mathcal{S}}^{\ast }$\textit{. Then}%
\begin{equation*}
\{\mathbf{u}_{x}\mid x\in \lbrack u]^{\ast }\}=\{\mathbf{u}_{x}\mid x\in
\lbrack u]^{1}\}\text{.}
\end{equation*}

\textit{In particular} 
\begin{equation*}
w_{u}=\underset{x\in \lbrack u]^{\ast }}{\vee }\mathbf{u}_{x}=\underset{x\in
\lbrack u]^{1}}{\vee }\mathbf{u}_{x}\text{.}
\end{equation*}

\textit{Proof}. Note that%
\begin{equation*}
\{\mathbf{u}_{x}\mid x\in \lbrack u]^{\ast }\}\subseteq \{\mathbf{u}_{x}\mid
x\in \lbrack u]^{1}\}\text{.}
\end{equation*}

Indeed, for $x\in \lbrack u]^{\ast }$ there exists $w_{x},y_{x}\in X$ such
that $x,y_{x}\in \mathcal{O}(w_{x})$ and $u(y_{x})=1$. Then, based on
Proposition 3.1, we have%
\begin{equation*}
\mathbf{u}_{x}=\underset{n\rightarrow \infty }{\lim }Z^{[n]}(\delta
_{y_{x}})=\mathbf{u}_{y_{x}}\text{,}
\end{equation*}%
so, as $u(y_{x})=1$, i.e. $y_{x}\in \lbrack u]^{1}$, we infer that $\mathbf{u%
}_{x}=\mathbf{u}_{y_{x}}\in \{\mathbf{u}_{x}\mid x\in \lbrack u]^{1}\}$.

As the inclusion 
\begin{equation*}
\{\mathbf{u}_{x}\mid x\in \lbrack u]^{1}\}\subseteq \{\mathbf{u}_{x}\mid
x\in \lbrack u]^{\ast }\}
\end{equation*}%
is obvious, the proof is completed. $\square $

\bigskip

\textbf{Proposition 3.7.} \textit{Let} $\mathcal{S}_{Z}=((X,d),(f_{i})_{i\in
I},(\rho _{i})_{i\in I})$ \textit{be an orbital fuzzy iterated function
system and }$u\in \mathcal{F}_{\mathcal{S}}^{\ast }$\textit{. Then }$w_{u}$ 
\textit{is upper semicontinuous, so }$w_{u}\in \mathcal{F}_{X}^{\ast }$.

\textit{Proof}. We have to show that%
\begin{equation}
\overline{\underset{n\rightarrow \infty }{\lim }}w_{u}(x_{n})\leq
w_{u}(x^{\ast })\text{,}  \tag{1}
\end{equation}
for every $(x_{n})_{n\in \mathbb{N}}\subseteq X$ and $x^{\ast }\in X$ such
that $\underset{n\rightarrow \infty }{\lim }x_{n}=x^{\ast }$.

It suffices to consider the case when $\overline{\underset{n\rightarrow
\infty }{\lim }}w_{u}(x_{n})\overset{not}{=}L>0$ since otherwise $(1)$ is
obvious.

Let us consider $(x_{n_{k}})_{k\in \mathbb{N}}$ such that $\underset{%
k\rightarrow \infty }{\lim }w_{u}(x_{n_{k}})=L$.

For the sake of simplicity, we denote $(x_{n_{k}})_{k\in \mathbb{N}}$ by $%
(z_{n})_{n\in \mathbb{N}}$.

So $\underset{n\rightarrow \infty }{\lim }z_{n}=x^{\ast }$ and $\underset{%
n\rightarrow \infty }{\lim }w_{u}(z_{n})=L>0$ and therefore we can suppose
that 
\begin{equation}
w_{u}(z_{n})\geq \frac{L}{2}\text{,}  \tag{2}
\end{equation}%
for every $n\in \mathbb{N}$.

Since%
\begin{equation*}
z_{n}\overset{(2)}{\in }[w_{u}]^{\frac{L}{2}}=[\underset{x\in \lbrack
u]^{\ast }}{\vee }\mathbf{u}_{x}]^{\frac{L}{2}}\overset{\text{Proposition 3.4%
}}{=}
\end{equation*}%
\begin{equation*}
=\underset{x\in \lbrack u]^{\ast }}{\cup }[\mathbf{u}_{x}]^{\frac{L}{2}}%
\overset{\text{Proposition 3.6}}{=}\underset{x\in \lbrack u]^{1}}{\cup }[%
\mathbf{u}_{x}]^{\frac{L}{2}}\text{,}
\end{equation*}%
for every $n\in \mathbb{N}$, we can consider a sequence $(y_{n})_{n\in 
\mathbb{N}}\subseteq \lbrack u]^{1}\subseteq $supp$u$ such that 
\begin{equation}
z_{n}\in \lbrack \mathbf{u}_{y_{n}}]^{\frac{L}{2}}\text{,}  \tag{3}
\end{equation}%
for every $n\in \mathbb{N}$.

As supp$u\in P_{cp}(X)$, there exists a subsequence $(y_{n_{k}})_{k\in 
\mathbb{N}}$ of $(y_{n})_{n\in \mathbb{N}}$ and $y\in X$ such that $\underset%
{k\rightarrow \infty }{\lim }y_{n_{k}}=y$.

The upper semicontinuity of $u$ implies%
\begin{equation}
u(y)\geq \overline{\underset{k\rightarrow \infty }{\lim }}u(y_{n_{k}})%
\overset{y_{n_{k}}\in \lbrack u]^{1}}{=}1>0\text{.}  \tag{4}
\end{equation}

\textbf{Claim}%
\begin{equation*}
w_{u}(z_{n_{k}})=\mathbf{u}_{y_{n_{k}}}(z_{n_{k}})\text{,}
\end{equation*}%
for every $k\in \mathbb{N}$.

\textit{Justification of the Claim}. Let us consider a fixed $k\in \mathbb{N}
$, but arbitrarily chosen. We have $w_{u}(z_{n_{k}})=\underset{x\in \lbrack
u]^{\ast }}{\vee }\mathbf{u}_{x}(z_{n_{k}})\overset{y_{n_{k}}\in \lbrack
u]^{1}}{\geq }\mathbf{u}_{y_{n_{k}}}(z_{n_{k}})\overset{(3)}{\geq }\frac{L}{2%
}>0$. Let us suppose, ad absurdum, that $w_{u}(z_{n_{k}})>\mathbf{u}%
_{y_{n_{k}}}(z_{n_{k}})$ and let us consider $\beta \in \mathbb{R}$ such
that 
\begin{equation}
w_{u}(z_{n_{k}})=\underset{x\in \lbrack u]^{\ast }}{\vee }\mathbf{u}%
_{x}(z_{n_{k}})>\beta >\mathbf{u}_{y_{n_{k}}}(z_{n_{k}})\text{.}  \tag{5}
\end{equation}%
Then there exists $s\in \lbrack u]^{\ast }$ such that 
\begin{equation}
\mathbf{u}_{s}(z_{n_{k}})>\beta >0\text{,}  \tag{6}
\end{equation}%
so, $z_{n_{k}}\in \lbrack \mathbf{u}_{s}]^{\ast }$ and, taking into account
Proposition 3.2, we infer that%
\begin{equation}
\mathbf{u}_{s}=\mathbf{u}_{z_{n_{k}}}\text{.}  \tag{7}
\end{equation}%
Since $y_{n_{k}}\in \lbrack u]^{1}$ and $\mathbf{u}_{y_{n_{k}}}(z_{n_{k}})%
\overset{(3)}{\geq }\frac{L}{2}>0$, Proposition 3.2 assures us that%
\begin{equation}
\mathbf{u}_{z_{n_{k}}}=\mathbf{u}_{y_{n_{k}}}\text{.}  \tag{8}
\end{equation}%
The contradiction $\beta \overset{(6)}{<}\mathbf{u}_{s}(z_{n_{k}})\overset{%
(7)\text{ \& }(8)}{=}\mathbf{u}_{y_{n_{k}}}(z_{n_{k}})\overset{(5)}{<}\beta $
ends the justification of the Claim.

Now let us consider a fixed $\varepsilon \in (0,\frac{L}{2})$, but
arbitrarily chosen.

As $\underset{k\rightarrow \infty }{\lim }w_{u}(z_{n_{k}})=L$, there exists $%
k_{\varepsilon }\in \mathbb{N}$ such that $w_{u}(z_{n_{k}})\overset{\text{%
Claim}}{=}\mathbf{u}_{y_{n_{k}}}(z_{n_{k}})\geq L-\varepsilon $, i.e.%
\begin{equation}
z_{n_{k}}\in \lbrack \mathbf{u}_{y_{n_{k}}}]^{L-\varepsilon }\text{,} 
\tag{9}
\end{equation}%
for every $k\in \mathbb{N}$, $k\geq k_{\varepsilon }$.

Since $\underset{k\rightarrow \infty }{\lim }y_{n_{k}}=y$, via Proposition
3.3, we deduce that $\underset{k\rightarrow \infty }{\lim }d_{\infty }(%
\mathbf{u}_{y_{n_{k}}},\mathbf{u}_{y})=0$ and therefore 
\begin{equation}
\underset{k\rightarrow \infty }{\lim }h([\mathbf{u}_{y_{n_{k}}}]^{L-%
\varepsilon },[\mathbf{u}_{y}]^{L-\varepsilon })=0\text{.}  \tag{10}
\end{equation}

As $\underset{k\rightarrow \infty }{\lim }z_{n_{k}}=x^{\ast }$, based on $%
(9) $, $(10)$ and Remark 2.1, we infer that $x^{\ast }\in \lbrack \mathbf{u}%
_{y}]^{L-\varepsilon }$, i.e. 
\begin{equation}
\mathbf{u}_{y}(x^{\ast })\geq L-\varepsilon >0\text{.}  \tag{11}
\end{equation}

A similar argument with the one used in the justification of the Claim
assures us that 
\begin{equation}
\mathbf{u}_{y}(x^{\ast })=w_{u}(x^{\ast })\text{.}  \tag{12}
\end{equation}

Indeed, we have $w_{u}(x^{\ast })=\underset{x\in \lbrack u]^{\ast }}{\vee }%
\mathbf{u}_{x}(x^{\ast })\overset{(4)}{\geq }\mathbf{u}_{y}(x^{\ast })$. Let
us suppose, ad absurdum, that $w_{u}(x^{\ast })>\mathbf{u}_{y}(x^{\ast })$
and let us consider $\beta \in \mathbb{R}$ such that $w_{u}(x^{\ast })=%
\underset{s\in \lbrack u]^{\ast }}{\vee }\mathbf{u}_{s}(x^{\ast })>\beta >%
\mathbf{u}_{y}(x^{\ast })$. Then there exists $s\in \lbrack u]^{\ast }$ such
that $\mathbf{u}_{s}(x^{\ast })>\beta >0$. We have $\mathbf{u}_{s}\overset{%
\text{Proposition 3.2}}{=}\mathbf{u}_{x^{\ast }}\overset{\text{Proposition
3.2}}{=}\mathbf{u}_{y}$ and therefore we get the contradiction $\beta <%
\mathbf{u}_{s}(x^{\ast })=\mathbf{u}_{y}(x^{\ast })<\beta $.

Thus, via $(11)$ and $(12)$ we infer that%
\begin{equation}
w_{u}(x^{\ast })\geq L-\varepsilon \text{.}  \tag{13}
\end{equation}

Relation $(13)$ implies the validity of $(1)$. $\square $

\bigskip

\textbf{Proposition 3.8.} \textit{Let} $\mathcal{S}_{Z}=((X,d),(f_{i})_{i\in
I},(\rho _{i})_{i\in I})$ \textit{be an orbital fuzzy iterated function
system and }$u\in \mathcal{F}_{\mathcal{S}}^{\ast }$\textit{. Then}%
\begin{equation*}
d_{\infty }(\mathbf{u}_{u},w_{u})=0\text{.}
\end{equation*}

\textit{Proof}. First let us note that, based on Remark 3.5, Lemma 3.4 from
[12] and Proposition 3.6, a repeated use of Remark 2.6 ensures us that%
\begin{equation}
\underset{x\in \lbrack u]^{\ast }}{\vee }Z^{[n]}(u^{x})=\underset{x\in
\lbrack u]^{\ast }}{\max }\text{ }Z^{[n]}(u^{x})\text{.}  \tag{1}
\end{equation}

We have%
\begin{equation*}
d_{\infty }(\mathbf{u}_{u},w_{u})\overset{\text{Remark 2.13, a)}}{\leq }%
d_{\infty }(\mathbf{u}_{u},Z^{[n]}(u))+d_{\infty }(\underset{x\in \lbrack
u]^{\ast }}{\vee }Z^{[n]}(u^{x}),\underset{x\in \lbrack u]^{\ast }}{\vee }%
\mathbf{u}_{x})\leq
\end{equation*}%
\begin{equation*}
\overset{\text{(1), Remark 2.7 \& Remark 3.5}}{\leq }d_{\infty }(\mathbf{u}%
_{u},Z^{[n]}(u))+\underset{x\in \lbrack u]^{\ast }}{\sup }d_{\infty
}(Z^{[n]}(u^{x}),\mathbf{u}_{x})\leq
\end{equation*}%
\begin{equation*}
\overset{\text{Remark 2.13, b)}}{\leq }d_{\infty }(\mathbf{u}%
_{u},Z^{[n]}(u))+\frac{C^{n}}{1-C}\underset{x\in \lbrack u]^{\ast }}{\sup }%
diam(F_{S}(\text{supp}u^{x})\cup \text{supp}u^{x})\leq
\end{equation*}%
\begin{equation}
\overset{\text{Remark 2.12}}{\leq }d_{\infty }(\mathbf{u}_{u},Z^{[n]}(u))+%
\frac{C^{n}}{1-C}diam(F_{S}(\text{supp}u)\cup \text{supp}u)\text{,}  \tag{2}
\end{equation}%
for every $n\in \mathbb{N}$.

As 
\begin{equation*}
\underset{n\rightarrow \infty }{\lim }d_{\infty }(\mathbf{u}_{u},Z^{[n]}(u))=%
\underset{n\rightarrow \infty }{\lim }\frac{C^{n}}{1-C}diam(F_{S}(suppu)\cup
suppu)=0\text{,}
\end{equation*}%
the conclusion follows by passing to limit, as $n$ goes to $\infty $, in $%
(2) $. $\square $

\bigskip

Now we state what can be viewed as the main theorem of this paper.

\bigskip

\textbf{Theorem 3.9.} \textit{Let} $\mathcal{S}_{Z}=((X,d),(f_{i})_{i\in
I},(\rho _{i})_{i\in I})$ \textit{be an orbital fuzzy iterated function
system and }$u\in \mathcal{F}_{\mathcal{S}}^{\ast }$\textit{. Then }%
\begin{equation*}
\mathbf{u}_{u}=\underset{x\in \lbrack u]^{\ast }}{\vee }\mathbf{u}_{x}=%
\underset{x\in \lbrack u]^{1}}{\vee }\mathbf{u}_{x}=\underset{x\in \lbrack
u]^{\ast }}{\max }\text{ }\mathbf{u}_{x}=\underset{x\in \lbrack u]^{1}}{\max 
}\text{ }\mathbf{u}_{x}\text{.}
\end{equation*}

\textit{Proof}. Since $\mathbf{u}_{u}\in \mathcal{F}_{X}^{\ast }$, $w_{u}%
\overset{\text{Proposition 3.7}}{\in }\mathcal{F}_{X}^{\ast }$, $d_{\infty }(%
\mathbf{u}_{u},w_{u})\overset{\text{Proposition 3.8}}{=}0$ and $d_{\infty }$
is a metric on $\mathcal{F}_{X}^{\ast }$ we conclude that $\mathbf{u}%
_{u}=w_{u}$, i.e., taking into account Proposition 3.6, we have%
\begin{equation*}
\mathbf{u}_{u}=\underset{x\in \lbrack u]^{\ast }}{\vee }\mathbf{u}_{x}=%
\underset{x\in \lbrack u]^{1}}{\vee }\mathbf{u}_{x}\text{. }\square
\end{equation*}%
Moreover, by virtue of Remark 3.5 and Proposition 3.6, we have%
\begin{equation*}
\mathbf{u}_{u}=\underset{x\in \lbrack u]^{\ast }}{\vee }\mathbf{u}_{x}=%
\underset{x\in \lbrack u]^{1}}{\vee }\mathbf{u}_{x}=\underset{x\in \lbrack
u]^{\ast }}{\max }\text{ }\mathbf{u}_{x}=\underset{x\in \lbrack u]^{1}}{\max 
}\text{ }\mathbf{u}_{x}\text{.}
\end{equation*}

\bigskip

\textbf{References}

\bigskip

[1] J. Andres, M. Rypka, Fuzzy fractals and hyperfractals, Fuzzy Sets and
Systems, \textbf{300} (2016), 40-56.

[2] M. Barnsly, S. Demko, Iterated function systems and the global
construction of fractals, Proc. Roy. Soc. London, \textbf{399} (1985),
243--275.

[3] M. Barnsley, Fractal Everywhere, Academic Press, Boston, 1988.

[4] C. Cabrelli, U. Molter, Density of fuzzy attractors: a step towards the
solution of the inverse problem for fractals and other sets, Probabilistic
and stochastic methods in analysis, with applications (ll Ciocco, 1991),
163-173. NATO Adv. Sci. Inst. Ser. C: Math. Phys. Sci., 372, Kluwer Acad.
Publ., Dordrecht. 1992.

[5] C. Cabrelli, B. Forte, U. Molter, E. Vrscay, Iterated fuzzy systems: a
new approach to the inverse problem for fractals and other sets, J. Math.
Anal. Appl., \textbf{171} (1992), 79-100.

[6] P. Diamond, P. Kloeden, Metric spaces of fuzzy sets, Theory and
applications, World Scientific Publishing Co., Inc., River Edge, NJ., 1994.

[7] B. Forte, M. Lo Sciavo, E. Vrscay, Continuity properties of attractors
for iterated fuzzy set systems, J. Austral. Math. Soc. Ser B., \textbf{36}
(1994), 175-193.

[8] J. Hutchinson, Fractals and self similarity, Indiana Univ. Math. J., 
\textbf{30} (1981), 713--747.

[9] R. Miculescu, A. Mihail, I. Savu, Iterated function systems consisting
of continuous functions satisfying Banach's orbital condition, An. Univ. de
Vest Timi\c{s}. Ser. Mat.-Inform., \textbf{56} (2018), 71-80.

[10] A. Mihail, I. Savu, Orbital $\varphi $-contractive iterated function
systems, Proceedings of Research World International Conference, Prague,
Czech Republic, 21-22 September 2020.

[11] A. Mihail, I. Savu, $\varphi $-contractive parent-child infinite IFSs
and orbital $\varphi $-contractive infinite IFSs, arXiv:2103.07551.

[12] A. Mihail, I. Savu, Orbital fuzzy iterated function systems,
arXiv:2112.15496.

[13] E. Oliveira, F. Strobin, Fuzzy attractors appearing from GIFZS, Fuzzy
Sets and Systems, \textbf{331} (2018), 131-156.

[14] I. Savu, New aspects concerning IFSs consisting of continuous functions
satisfying Banach's orbital condition, J. Fixed Point Theory Appl., \textbf{%
21} (2019), 62.

[15] R. Uthayakumar and D. Easwaramoorthy, Hutchinson-Barnsley operator in
fuzzy metric spaces, Int. J. Math. Comput. Sci., World Acad. Sci., Eng.
Technol., \textbf{5} (2011), 1418-1422.

[16] L. Zadeh, Fuzzy sets, Information and Control, \textbf{8} (1965),
338-353.

\bigskip

{\small Radu MICULESCU}

{\small Faculty of Mathematics and Computer Science}

{\small Transilvania University of Bra\c{s}ov}

{\small Iuliu Maniu Street, nr. 50, 500091, Bra\c{s}ov, Romania}

{\small E-mail: radu.miculescu@unitbv.ro}

\bigskip

{\small Alexandru MIHAIL}

{\small Faculty of Mathematics and Computer Science}

{\small Bucharest University, Romania}

{\small Str. Academiei 14, 010014, Bucharest}

{\small E-mail: mihail\_alex@yahoo.com}

\bigskip

{\small Irina SAVU}

{\small Faculty of Applied Sciences}

{\small Politehnica University of Bucharest}

{\small Splaiul Independen\c{t}ei 313, Bucharest, Romania}

{\small E-mail: maria.irina.savu@upb.ro}

\end{document}